\documentclass{article}
\usepackage{graphicx}
\usepackage{amsmath}
\usepackage{amsfonts}
\usepackage{amsthm}

\newtheorem{theorem}{Theorem}[section]
\newtheorem{lemma}[theorem]{Lemma}

\begin{document}

\title{Unique Tournaments and Radar Tracking}

\author{Tanya Khovanova}
\date{December 4, 2007}
\maketitle

\begin{abstract}
The sequence counting the number of unique tournaments with $n$ people is the same as the sequence counting non-tracking binary strings corresponding to $n-2$ radar observations with the tracking rule ``3 out of 5 with loss 2.'' This fact allows us to build a bijection between unique tournaments and non-tracking binary strings.
\end{abstract}


\section{Introduction}

This paper is about an integer sequence that appears in two very different places. The sequence in question is sequence 
A000570 in the Online Encyclopedia of Integer Sequences (OEIS), where it is described as: ``Number of tournaments on n nodes determined by their score vectors.'' It is interesting to note that I generated this sequence while playing with rules for radar tracking. The material in this paper loosely follows  the order in which I discovered the facts about this sequence.

\section{Radar Tracking with Loss}

Radars or sonar detectors detect and track targets in a series of observations. The result of an observation can be described as detection or no detection. There are rules for a series of detections to produce a track on a target. One of the most commonly used rules is called ``3 out of 5 with 2 loss'' \cite{Brookner}. The rule means that to get a track, the radar must have at least 3 detections in a series of 5 or fewer observations and there shouldn't be two non-detections in a row between two consecutive detections. In other words, if you have two consecutive non-detections you have lost your track and you must start again.

If we represent a detection as 1 and a non-detection as 0, the result of series of $k$ observations is a binary string of length $k$. For every binary string we can say whether the detections and non-detections described by this string establish a track or not. If we count the number of track-producing strings of length $k$, we get an integer sequence that we call a tracking sequence. Correspondingly, we can count the number of non-track producing strings of length $k$ and get a non-tracking sequence. Obviously, we can define tracking and non-tracking sequences for any tracking rule. In this paper we only discuss the particular rule ``3 out of 5 with loss 2.''

Let me introduce some notations. I denote by $Tr(k)$ the tracking sequence --- the number of binary strings of length $k$ that produce a track.  The non-tracking sequence $NTr(k)$ is the number of binary strings of length $k$ that do not allow tracking. 

Obviously, $Tr(k) + NTr(k) = 2^k$. That means to calculate these two sequences, it is enough to calculate one of them. Non-tracking sequences are usually easier to calculate. So, let's calculate $NTr(k)$.

Here is the list of non-tracking strings for $k \le 5$. 
\begin{enumerate}
\item For $k = 1$ any string is a non-tracking string. The list of non-tracking strings is: 0, 1. The total is 2.
\item The same thing happens for $k=2$. The list of non-tracking strings is:  00, 01, 10, 11. The total is 4.
\item For $k = 3$, the only tracking string is the string with 3 ones. Hence the non-tracking strings are: 000, 001, 010, 011, 100, 101, 110. The total is 7.
\item For $k = 4$ a non-tracking string is any string of length 4 that doesn't contain 3 ones: 0000, 0001, 0010, 0011, 0100, 0101, 0110, 1000, 1001, 1010, 1100. The total is 11.
\item For $k = 5$ a non-tracking string is any string of length 5 that doesn't contain 3 ones, plus the strings that contain 3 ones but have two zeroes in between. Namely: 00000, 00001, 00010, 00011, 00100, 00101, 00110, 01000, 01001, 01010, 01100, 10000, 10001, 10010, 10011, 10100, 11000, 11001. The total is 18.
\end{enumerate}
We see that $NTr(k)$ starts as: 2, 4, 7, 11, 18, $\ldots$. I will not bore you by describing the Java program I wrote to calculate this sequence. The result was: 2, 4, 7, 11, 18, 31, 53, 89, 149, 251, 424, 715, 1204, 2028, $\ldots$.

I plugged the numbers into the Online Encyclopedia of Integer Sequences (OEIS) \cite{OEIS} and I got: A000570 --- Number of tournaments on n nodes determined by their score vectors.

\section{Unique Tournaments}

In mathematics a tournament describes the results of a bunch of people who gather together to play a game of tennis in a round-robin or all-play-all style. You can replace tennis by any other game where every two people play with each other once and draws are not allowed. So, a mathematical tournament is a competition in which there is a winner for every pair of competitors. Or to put it more formally, a tournament is a complete oriented graph.

The score vector is the list of total wins for every participant, where the order of participants does not matter. In other words, the score vector can be defined as a non-decreasing array of outdegrees of the vertices of the tournament graph. 
A tournament is called unique if it is uniquely defined by its score vector.

I will denote by $UT(k)$ the number of unique tournaments with $k$ people (vertices). Let us take a look at the unique tournaments sequence:

1, 1, 2, 4, 7, 11, 18, 31, 53, 89, 149, 251, 424, 715, 1204, 2028, 3418, $\ldots$.

\begin{figure}
  \centering
  \includegraphics[scale=0.5]{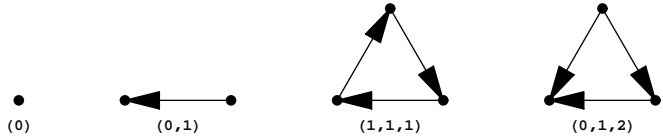}
  \caption{Small Tournaments.}\label{small}
\end{figure}

Here the first element has index 1. There is only one tournament with only 1 participant: no games were played and the score vector is $\{0\}$. With two participants there exists only one tournament: one person wins and the other person loses. The score vector of this tournament is $\{0, 1\}$. With 3 participants there can be 2 different tournaments. The first tournament with a score vector $\{0, 1, 2\}$ has one champion who won every game and one loser who lost every game.  The second tournament with a score vector $\{1, 1, 1\}$ doesn't have a champion. See Figure~\ref{small} for the graphs of very small tournaments of up to three people.

\begin{figure}
  \centering
  \includegraphics[scale=0.5]{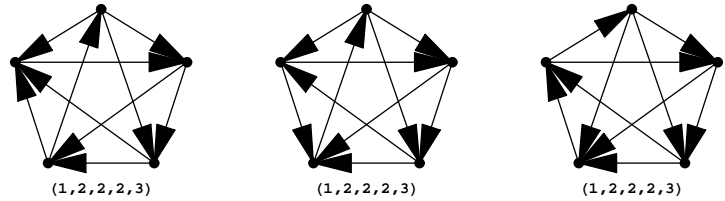}
  \caption{Non-Unique Tournaments.}\label{nonunique}
\end{figure}

Tournaments up to 4 people are always unique. Non-unique tournaments appear when there are 5 competitors or  more. For tournaments with 5 people there are two score vectors, $\{1,2,2,2,3\}$ and $\{1,1,2,3,3\}$, that define several non-isomorphic tournaments. Figure~\ref{nonunique} shows 3 non-isomorphic tournaments defined by the score vector $\{1,2,2,2,3\}$.

The indexing of the unique tournament sequence is slightly different from the tracking sequence. Namely, $NTr(k) = UT(k+2)$. We will return to this shift later. 

If we check the entry in the OEIS about the unique tournaments sequence we can find a recurrence relation for it:
\begin{equation}\label{E:rec}
a(n) = a(n-1) + a(n-3) + a(n-4) + a(n-5).
\end{equation}

I have shown that the non-tracking sequences and the unique tournament sequences have the same start. To prove that these sequences are truly the same, it is enough to prove the recurrence relation \ref{E:rec} for the non-tracking sequence.

\section{Recurrence Relation}

\begin{lemma}\label{L:rec} The non-tracking sequence $NTr(k)$ satisfies the recurrence: 
$$NTr(k) = NTr(n-1) + NTr(n-3) + NTr(n-4) + NTr(n-5).$$
\end{lemma}

\begin{proof}
Let us assume that $k$ is more than 5 and count the number of non-tracking strings of length $k$. Every non-tracking string of length $k$ can be one of the following:
ends with 0, ends with 01 or ends with 11. The number of non-tracking strings of length $k$ that end with 0 is exactly the same as the number of
all non-tracking strings of length $k-1$. If a non-tracking string ends in 11 it has to have
two zeroes before 11, that is, it has to actually end in 0011. The number
of non-tracking strings of length $k$ that end in 0011 is obviously $NTr(k-4)$. The non-tracking strings that end in 01 can be of two types: end in 001 and end in 101. If a string ends in 001, that means the track is just lost and before that there can be any non-tracking string. That is, the number of non-tracking strings ending in 001 is $NTr(k-3)$. If the non-tracking string ends in 101 it has to have two zeroes before that. That means it really ends in  00101. Two zeroes signify the loss of the track and the string can have anything before that. That means the number of non-tracking strings ending in 101 is $NTr(k-5)$. By summing over all the endings we get the proof.
\end{proof}

We proved that $NTr(k) = UT(k+2)$. The sequences are the same, but can we match a specific tournament to a non-tracking binary string? 

\section{Matching Tournaments and Binary Strings}

\begin{figure}
  \centering
  \includegraphics[scale=0.5]{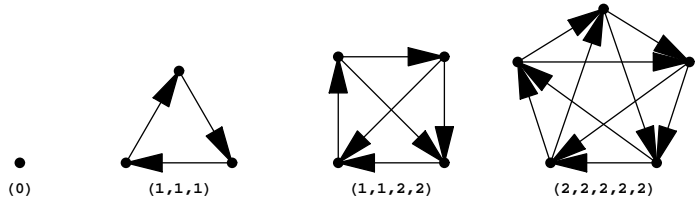}
  \caption{Basic Tournaments.}\label{basic}
\end{figure}

Prasad Tetali \cite{tetali} shows us that there are four basic unique tournaments. Their score vectors are: $\{0\}$, $\{1, 1, 1\}$, $\{1, 1, 2, 2\}$ and $\{2, 2, 2, 2, 2\}$. See Figure~\ref{basic} for the pictures of unique tournaments.
We can build other tournaments from basic tournaments using the composition operation. Suppose we have two tournaments $A$ and $B$ with $m$ and $n$ nodes correspondingly. We can create a new tournament $A + B$ on $m+n$ nodes in the following manner: The first group of $m$ people played with each other according to the tournament $A$. The second group of $n$ people played with each other according to the tournament $B$. And every person from the second group won every match with a person from the first group. According to \cite{tetali} any unique tournament can be decomposed into basic unique tournaments.

As we are trying to draw a parallel between tournaments and non-tracking strings, we have candidates for building blocks in non-tracking strings. Namely, there are 4 basic strings that we add to an existing non-tracking string to create a new one: 0, 001, 0011, 00101. If these strings correspond to basic unique tournaments then we can suppose that one character in a binary string should correspond to one person (vertex) in a tournament. But something is not quite right here. Actually two things are off:
\begin{enumerate}
\item The index is off --- tournaments with $k+2$ people correspond to non-tracking strings of length $k$.
\item The initial non-tracking sequences can't be decomposed into basic strings. For example, a non-tracking sequence ``1'' of length 1 is not a concatenation of basic non-tracking sequences.
\end{enumerate}

Now is time for some magic. Let us create a new sequence. This sequence will be denoted as $ILNTr(k)$ -- the sequence of initial-loss non-tracking strings of length $k$. Initial-loss non-tracking strings are non-tracking strings that can't have ones in the first and second place. In other words, the initial radar observations start by losing the track. So, we have only one initial-loss non-tracking string of length 1: 0. We also have only one such string of length 2: 00. The number of initial-loss non-tracking strings of length $k \ge 2$ is equal to the number of non-tracking strings of length $k-2$. Indeed, we can get an initial-loss non-tracking string of length $k \ge 2$ by adding two zeroes to any non-tracking string of length $k-2$.

\begin{lemma} Any initial-loss non-tracking binary string is a concatenation of strings of 4 basic types: 0, 001, 0011, 00101. 
\end{lemma}

\begin{proof} An initial-loss non-tracking string is a valid non-tracking string. By the proof of lemma \ref{L:rec} if our initial-loss non-tracking string has a length more than 5, it has to end in a string of one of the basic types. Thus, we can remove all the strings of the basic type from the right of our string until this string has a length not more than 5. For the strings of small length we can manually check each one. Before starting manual labor let me simplify my life a little bit more. Namely, zeroes at the end of the string are concatenations of the first basic type: 0. So, I can only manually check initial-loss non-tracking strings of length not more than 5 ending in 1. Here is the table providing the decomposition:

\begin{center}
\begin{tabular}[c]{|c|c|c|}
\hline
non-tracking & IL non-tracking & decomposition \\
\hline
1    & 001   & 001  \\
01   & 0001  & 0 + 001 \\
11   & 0011  & 0011    \\
001  & 00001 & 0 + 0 + 001 \\
011  & 00011 & 0 + 0011 \\
101  & 00101 & 00101 \\
\hline
\end{tabular}
\end{center}

Thus we established a one-to-one correspondence between unique tournaments and initial-loss non-tracking binary strings. 
\end{proof}

Now we can look at some properties of unique tournaments and try to translate these into properties for binary strings. For example, there is a natural duality for tournaments. Can we extend the duality to initial-loss non-tracking strings?

\section{Duality}

There is a natural duality of tournaments. Given a tournament the dual tournament is the tournament on the same number of nodes with directions of all edges reversed. In other words, for every pair of people the person who wins in the initial tournament loses in the dual tournament.

Obviously if a tournament is unique then its dual is a unique tournament too. It would be interesting to transfer this duality to initial-loss non-tracking binary strings.

Let us first see how the duality acts on a score vector. Obviously, if a person won $x$ games in an $n$ persons tournament, he would win $n-1-x$ games in a dual tournament. That means that the score vector of the dual tournament can be calculated by replacing every score $x$ in the original tournament by $n-1-x$. 

\begin{lemma}
All basic unique tournaments are self-dual.
\end{lemma}

\begin{proof}
We can check the lemma directly by reversing the arrows on the tournaments graphs of basic tournaments in Figure~\ref{basic}. Another way to prove it is to check that the score vectors of the duals to basic tournaments do not change. The proof follows from the fact that the tournament is defined by its score.
\end{proof}

Now I would like to discuss how the duality affects the decomposition of tournaments. The dual tournament to the composition $A+B$ is the composition of the tournament dual to $B$ with the tournament dual to $A$.

From this it is easy to see that the dual to a unique tournament can be constructed by decomposition of this tournament into basic tournaments then composing them back together in the reverse order.

From here we see how the duality works for initial-loss non-tracking strings. Take one such string, decompose it into a concatenation of basic strings, then compose them back in the reverse order.

\section{Score Vector}

Can we look at an initial-loss non-tracking binary string and recover the score vector of the corresponding tournament directly from the string? The strings that represent a unique tournament can be decomposed into a concatenation of basic strings. Basic strings correspond to basic tournaments. If we know the decomposition of a tournament into basic tournaments we can calculate the score vector. For example, suppose we decompose a tournament into two basic tournaments $A$ and $B$. The decomposition means that all participants of the tournament $B$ won their games with any participant of the tournament $A$. Hence, the score vector for participants from $A$ is the same as the score vector for the tournament A. The score vector for the participants from $B$ are the score vector of $B$ plus the number of participants in $A$.

Now I will describe the score vector given an initial-loss non-tracking binary string of length $k$. Given index $i$, I would like to find the score of an $i$-th person/node. Let us denote by $n$ the index of the start of the first occurrence of two zeroes in our string that is bigger than $i$. If $i$ is close to the end and there are no 2 zeroes after $i$, then $n = k+1$. Let us denote by $m$ the index of the start of the latest occurrence of two zeroes before or including $i$. It is easy to see that the substring starting from $m$ and ending with $n-1$ is a basic substring. The score to the i-th person is $i$ plus an adjustment coefficient which depends on the place of the i-th character in the basic substring. The adjustment number is given by the following table:

\begin{center}
\begin{tabular}[c]{|c|c|}
\hline
basic string & vector of score adjustments \\
\hline
0   & \{0\} \\
001  &  \{1, 0, -1\} \\
0011  & \{1, 0, 0, -1\} \\
00101 & \{2, 1, 0, -1, -2\} \\
\hline
\end{tabular}
\end{center}

In particular, if the string corresponding to the tournament ends with 0, we know that there is a champion who won all the matches. On the other hand if the binary string starts with three zeroes we know that there is an absolute loser.

\section{Future Research}

There are several questions for future research:

\begin{enumerate}
\item What happens if we generalize the sequence to other tracking rules?
\item Can the matching of binary strings to unique tournaments be extended to include other tournaments or other score vectors?
\item Can some properties of binary strings be easily described in terms of the corresponding unique tournaments, for example the number of ones in an initial-loss non-tracking binary string?
\item Can some properties of unique tournaments be easily described in terms of the corresponding binary strings?
\end{enumerate}

\section{Acknowledgements}

I would like to thank OEIS (Online Encyclopedia of Integer Sequences) for being an extremely helpful resource in studying integer sequences. I am thankful to Alexey Radul, Jane Sherwin and Sue Katz for helping me with English for this paper.

\end{document}